\documentclass[pdflatex,sn-mathphys-num]{sn-jnl}

\usepackage{amsmath,amssymb,amsfonts}
\usepackage{amsthm}
\usepackage{booktabs}
\usepackage{xcolor}
\usepackage{url}
\usepackage{tikz}

\theoremstyle{thmstyleone}
\newtheorem{theorem}{Theorem}[section]
\newtheorem{proposition}[theorem]{Proposition}
\newtheorem{lemma}[theorem]{Lemma}
\newtheorem{corollary}[theorem]{Corollary}

\newtheorem{maintheorem}{Theorem}

\theoremstyle{thmstyletwo}
\newtheorem{example}[theorem]{Example}

\theoremstyle{thmstylethree}
\newtheorem{definition}[theorem]{Definition}

\newcommand{\Ap}{\operatorname{Ap}}
\newcommand{\Fr}{\operatorname{Fr}}
\newcommand{\Tot}{\operatorname{Tot}}
\newcommand{\reg}{\operatorname{reg}}
\newcommand{\M}{\operatorname{M}}
\newcommand{\N}{\mathbb{N}}
\newcommand{\Z}{\mathbb{Z}}

\begin{document}

\title[Undershoot of $h$-fold sumsets]{Undershoot of $h$-fold sumsets: Ap\'ery spectra, similar ordering and the three-distance theorem}

\author*[1]{\fnm{Prashant} \sur{Chaudhary}}\email{prashant@prashantchaudhary.com}
\author[2]{\fnm{Kapil} \sur{Khandelwal}}\email{kapil.khandelwal@nd.edu}

\affil*[1]{\orgname{Independent Researcher}, \orgaddress{\city{San Francisco Bay Area}, \state{CA}, \country{USA}}}
\affil[2]{\orgdiv{Department of Civil and Environmental Engineering and Earth Sciences}, \orgname{University of Notre Dame}, \orgaddress{\city{Notre Dame}, \state{IN}, \country{USA}}}

\abstract{For a finite set $A_0\subset\mathbb N$ with $0\in A_0$, $\max A_0=d$ and $\gcd(A_0)=1$, and a target $T$, the undershoot at level $h$ is the distance from $T$ down to the largest element of the $h$-fold sumset $hA_0$ not above $T$. Its total over $T\in[0,hd]$ is bounded in $h$: for large $h$ it equals the sum of the masses of two numerical semigroups attached to $A_0$, where the mass adds $L(L+1)/2$ over the maximal runs of $L$ gaps. We describe how the total reaches this limit: it can rise to it, fall to it, or climb past it and return. When the monomial curve of $A_0$ is arithmetically Cohen--Macaulay, two Ap\'ery sets give every undershoot exactly, and partial sums of a fixed spectrum bound the total from below, with equality at every level exactly when two threshold sequences are similarly ordered; the total then rises concavely to its limit, reached exactly at the sumsets regularity. We characterize the similarly ordered four-element sets. For $\{0,a,b,2b\}$ with $0<a<b$ coprime, the three-distance theorem gives the total in explicit quadratic pieces. A non-Cohen--Macaulay set exceeds its limit at every failing level $h$ with $hd$ at least the sum of the two conductors; every such four-element set does so, and a five-element example never does. For arbitrary sets, finite corner data give the exact undershoot over any window.}



\maketitle

\section{Introduction}\label{sec:intro}

The sums of at most $h$ elements of a finite set of integers need not fill the range below $h$ times its largest element. This paper measures how far they fall short: for each target, the distance down to the nearest sum not above it, added over all targets, as a function of $h$. Let $A_0=\{0=a_0<a_1<\dots<a_k=d\}$ be a set of integers with $\gcd(A_0)=1$, and let its $h$-fold sumset $hA_0\subseteq[0,hd]$ be the set of sums of at most $h$ elements of $A=A_0\setminus\{0\}$. For a target $T\ge0$ the \emph{undershoot}
\[
D_h(T)=T-\max\{n\in hA_0:n\le T\}
\]
is the distance from $T$ down to the largest sum not above it, and $\Tot(h)=\sum_{T=0}^{hd}D_h(T)$ is its total over the targets in range. The question came from filling advertising breaks: if the elements of $A$ are the available ad durations and a break of length $T$ takes at most $h$ ads, then $D_h(T)$ is the time left unfilled when the break is filled as fully as possible. The range of targets grows with $h$, yet the total stays bounded and is eventually constant. How it reaches its final value depends on the set. For $\{0,5,7,14\}$ the totals at $h=1,2,\dots$ are $32,48,56,61,64,65,65,\dots$: they rise by smaller steps each time and then stop. For $\{0,2,8,11\}$ they are $19,19,18,18,\dots$: they start above the final value and come down to it. For $\{0,4,9,13\}$ they are $22,36,45,48,49,46,44,40,40,\dots$: they climb past the final value $40$ and return to it only at $h=8$. This paper explains the three behaviors and describes the whole sequence $\Tot(h)$.

\paragraph{The limit.} For large $h$ the set $hA_0$ is an interval with fixed sets of exceptions at both ends \cite{Nathanson1972,GranvilleShakan2020,GranvilleWalker2021,Lev2022}, governed by the two numerical semigroups $S_1=\langle A\rangle$ and $S_2=\langle d-A_0\rangle$. Through the projective monomial curve of $A_0$, the level from which this holds is a regularity \cite{Elias2022,GimenezGonzalezSanchez2023}, and the sizes $|hA_0|$ have recently been computed, for special classes of sets \cite{Nathanson2025} and, from the strict records of the least number of parts in each residue class, for every finite set \cite[Thm.~1]{ChenYang2026}. These results count the sums. The undershoot records how far each target is from one. The end structure gives its limit: for every $A_0$ and all large $h$, $\Tot(h)=\M(S_1)+\M(S_2)$ (Propositions~\ref{prop:failing-interval} and~\ref{prop:overshoot}), where the mass $\M(S)$ of a numerical semigroup $S$ is the sum of $L(L+1)/2$ over its deserts, the maximal runs of gaps, of lengths $L$. For the three sets above the limits are $65$, $18$ and $40$. The end structure does not show how the totals approach the limit, and that is the subject of this paper.

\paragraph{Main results.} Let $c_1$ and $c_2$ be the conductors of $S_1$ and $S_2$, and let
\[
I_h=\{n\in[0,hd]:n\in S_1,\ hd-n\in S_2\}.
\]
Always $hA_0\subseteq I_h$. We call $A_0$ \emph{Cohen--Macaulay} (CM) when $hA_0=I_h$ for every $h$; by Cavaliere and Niesi this is the arithmetic Cohen--Macaulay property of the monomial curve (Theorem~\ref{thm:cm}). For a CM set and each residue $s$ modulo $d$ there are integers $u_s,v_s\ge0$, read off the Ap\'ery sets of $S_1$ and $S_2$, such that $kd-s\in hA_0$ exactly when $u_s\le k\le h-v_s$ (Lemma~\ref{lem:residue}). We call the pairs $(u_s,v_s)$ \emph{similarly ordered} if $u_s>u_t$ implies $v_s\ge v_t$, the hypothesis of Chebyshev's sum inequality. The sumsets regularity $\sigma(A)$ of Gimenez and Gonz\'alez-S\'anchez is the least level from which $hA_0$ keeps its final shape (Section~\ref{sec:prelim}).

\begin{maintheorem}[Theorem~\ref{thm:spectrum}, Proposition~\ref{prop:threshold}]\label{thm:A}
Let $A_0$ be CM with $d\ge2$. There are integers $\lambda_1\ge\lambda_2\ge\dots\ge0$, finitely many of them positive, with $\sum_\ell\lambda_\ell=\M(S_1)+\M(S_2)$, such that
\[
\Tot(h)\ge\lambda_1+\dots+\lambda_h\qquad(h\ge1),
\]
with equality for every $h$ if and only if the pairs $(u_s,v_s)$ are similarly ordered. In that case $\Tot$ is nondecreasing and concave, $\Tot(h)=\M(S_1)+\M(S_2)$ exactly when $h\ge\sigma(A)$, and the Castelnuovo--Mumford regularity of the coordinate ring of the curve is $\sigma(A)+1$.
\end{maintheorem}

Proposition~\ref{prop:threshold} proves the threshold and the regularity statement at every modulus $m\ge2$, a divisor of $d$ lying in $S_1\cap S_2$, at which the pairs read modulo $m$ are similarly ordered. The set $\{0,5,7,14\}$ is similarly ordered with $\lambda=(32,16,8,5,3,1)$, and its totals are the partial sums. The set $\{0,2,8,11\}$ is not: its bound $10,15,18$ is strict at $h=1$ and $h=2$, and the levels at which the bound is strict are given in general by the inversions of the pairs (Proposition~\ref{prop:strict-levels}). For four elements similar ordering has a finite geometric test. The pairs of a CM set $\{0,a_1,a_2,d\}$ come from its Ap\'ery diagram, a set of $d$ points $(x,y)\in\N^2$ such that the values $a_1x+a_2y$ are the least elements of $S_1$ in the $d$ residue classes (Lemma~\ref{lem:diagram}).

\begin{maintheorem}[Theorem~\ref{thm:similar}, Proposition~\ref{prop:shape}]\label{thm:B}
For a CM set $\{0,a_1,a_2,d\}$ the Ap\'ery diagram is a rectangle or an L-shape, and the pairs are similarly ordered if and only if on each layer $x+y=t$ of the diagram the values $a_1x+a_2y$ lie in one interval $((q-1)d,qd]$ with $q\in\Z$. An L-shape whose notch has both sides of length at least $2$ is never similarly ordered.
\end{maintheorem}

The family $\{0,a,b,2b\}$ is CM. Its pairs need not be similarly ordered at $d=2b$, but they always are when the residues are taken modulo $b$, and its spectrum is then a statistic of an arithmetic progression on a circle. For $1\le j\le b$ let $\phi(j)$ be the sum of $\binom g2$ over the gaps $g$ of the points $0,a,\dots,(j-1)a$ on $\Z/b\Z$.

\begin{maintheorem}[Theorems~\ref{thm:family-spectrum} and~\ref{thm:three-distance}, Corollaries~\ref{cor:family-mass} and~\ref{cor:quadratic-pieces}]\label{thm:C}
Let $0<a<b$ be coprime and $A_0=\{0,a,b,2b\}$. Then
\[
\Tot(h)=\binom b2+2\sum_{j=2}^{h}\phi(j)+\phi(h+1)\qquad(1\le h\le b-1),
\]
and $\Tot(h)=2\sum_{j=1}^b\phi(j)-\binom b2$ for $h\ge b-1$. The gaps of $\{0,a,\dots,(j-1)a\}$ take at most three lengths, given by the continued fraction of $a/b$ through the three-distance theorem, and $\Tot$ consists of explicit quadratic pieces, one for each block of that continued fraction.
\end{maintheorem}

For $\{0,5,7,14\}$, $\phi(1),\dots,\phi(6)=21,11,5,3,2,1$, which gives $32,48,56,61,64,65$. A set that is not CM has \emph{failing levels}, the $h$ with $hA_0\neq I_h$. They are finitely many, and $\tau$ denotes the level after the last one.

\begin{maintheorem}[Propositions~\ref{prop:failing-interval} and~\ref{prop:overshoot}, Corollary~\ref{cor:four-overshoot}, Example~\ref{ex:five-element}]\label{thm:D}
Let $A_0$ be arbitrary and $hd\ge c_1+c_2$. Then $\Tot(h)\ge\M(S_1)+\M(S_2)$, with equality if and only if $h$ is not a failing level. Every four-element set that is not CM has $\Tot(\tau-1)>\M(S_1)+\M(S_2)$. The five-element set $\{0,5,9,11,20\}$ is not CM, and $\Tot(h)\le\M(S_1)+\M(S_2)$ for every $h$.
\end{maintheorem}

The set $\{0,4,9,13\}$ has $c_1+c_2=48$ and failing levels $4$ to $7$, so by Theorem~\ref{thm:D} its total exceeds $40$ at each of them and equals $40$ from $h=8$ on. The four-element statement uses the companion paper \cite[Theorem~A]{CompanionRegularity}, which shows that $(\tau-1)d\ge c_1+c_2+d$ for these sets. For arbitrary sets the undershoot differs from that of $I_h$ by an explicit correction (Proposition~\ref{prop:correction}), and finitely many corners per residue class give the exact undershoot over any window of targets (Theorem~\ref{thm:hinge}).

\paragraph{Method.} Everything is read off the residue pairs $(u_s,v_s)$, which locate each residue class in $hA_0$ block by block. For a CM set the integers of $[0,hd]$ missing from $hA_0$ fall in blocks between consecutive multiples of a modulus, and the undershoot of a target is the number of consecutive missing integers ending at it (Proposition~\ref{prop:counting}). The missing positions of a block are the union of a set coming from $S_1$ and a set coming from $S_2$. Taking the larger of the two in every block gives the spectral bound, and the bound is attained in a block exactly when the two sets are nested. At the modulus $d$ they are nested in every block at every level exactly under similar ordering. For four elements the pairs come from the Ap\'ery diagram, and the test becomes a test on its layers. For $\{0,a,b,2b\}$ the missing positions at the modulus $b$ are the complements of the progressions $0,-a,\dots,-(j-1)a$ on $\Z/b\Z$, so the spectrum is a sum over circular gaps and the three-distance theorem evaluates it. For a set that is not CM we compare $hA_0$ with $I_h$: once $hd\ge c_1+c_2$ the total of $I_h$ is the mass, and every element of $I_h$ missing from $hA_0$ raises the total. The corrections and window formulas come from the minimal corners of the levels at which each class is present.

\paragraph{Related work.} The three-distance theorem is due to S\'os \cite{Sos1958}, Sur\'anyi \cite{Suranyi1958} and \'Swierczkowski \cite{Swierczkowski1958}; we use it in the explicit form of Alessandri and Berth\'e \cite{AlessandriBerthe1998}. The homogenized Ap\'ery set and the Cohen--Macaulay criterion are those of Cavaliere and Niesi \cite{CavaliereNiesi1983} and Gimenez and Gonz\'alez-S\'anchez \cite{GimenezGonzalezSanchez2023}. Lemma~\ref{lem:residue} is their coordinate description read on one level, and Proposition~\ref{prop:counting} follows from it directly. Proposition~\ref{prop:failing-interval} is \cite[Lemma~2.5]{CompanionRegularity} in other coordinates, and the L-shapes of Section~\ref{sec:similar} are those of three-generated numerical semigroups \cite{AguiloMarijuan2014}. What this paper adds begins with the spectrum. Hoa and Tien \cite{HoaTien2026} bound the regularity of the monomial curve through Ap\'ery sets and Frobenius numbers, with criteria for the Cohen--Macaulay property and an application to the structure of sumsets. Similar ordering is the condition of Chebyshev's sum inequality \cite{HardyLittlewoodPolya1952}. To our knowledge the undershoot has not been studied as a functional of $hA_0$. The main question left open is a description of the similarly ordered four-element sets by a short list of congruence or continued-fraction conditions on $(a_1,a_2,d)$ (Section~\ref{sec:open}).

\paragraph{Organization.} Section~\ref{sec:prelim} fixes notation. Sections~\ref{sec:sumsets} to~\ref{sec:noncm} prove Theorems~\ref{thm:A} to~\ref{thm:D} in order, and Section~\ref{sec:corrections} gives the corrections and window formulas. Section~\ref{sec:open} lists open problems, and Section~\ref{sec:computation} describes the computations, all exact and reproducible from the accompanying code \cite{code2026}.

\section{Preliminaries}\label{sec:prelim}

\paragraph{Numerical semigroups.} Here $\N=\{0,1,2,\dots\}$, and an interval $[a,b]$ means the set of integers $n$ with $a\le n\le b$ (likewise for half-open intervals), except on the circle $\mathbb R/\Z$ in Section~\ref{sec:family}. A numerical semigroup is a submonoid $S$ of $(\N,+)$ with finite complement. Its elements not in $S$ are its gaps, its largest gap is the Frobenius number $\Fr(S)$, and $c(S)=\Fr(S)+1$ is its conductor; for $S=\N$ we put $\Fr(S)=-1$ and $c(S)=0$. For a positive element $m$ of $S$, the Ap\'ery set of $S$ with respect to $m$ is $\Ap(S;m)=\{w_r:r\in\Z/m\Z\}$, where $w_r$ is the least element of $S$ congruent to $r$ modulo $m$; an integer $n\ge0$ lies in $S$ if and only if $n\ge w_{n\bmod m}$, and $\Fr(S)=\max_r w_r-m$ \cite{RosalesGarciaSanchez2009}.

\begin{definition}[Mass]\label{def:mass}
The \emph{mass} of a numerical semigroup $S$ is
\[
\M(S)=\sum_{x\notin S}\bigl(x-\max\{s\in S: s<x\}\bigr),
\]
the sum over the gaps of $S$ of the distance down to the preceding element of $S$.
\end{definition}

A desert, that is, a maximal run of consecutive gaps \cite{NumericalSgps} (for two generators see \cite[Sect.~6]{Moree2014}), of length $L$ contributes $1+2+\dots+L=L(L+1)/2$, and so does the sum of the distances up to the next element. For a finite set $X\subset\Z$ write $f(X)$ for the number of pairs $s\le t$ with $[s,t]\subseteq X$; if $X$ is a union of maximal runs of lengths $L_i$, then $f(X)=\sum_i L_i(L_i+1)/2$. Thus $\M(S)=f(\N\setminus S)$.

\paragraph{Circular gaps.} We use one statistic of arithmetic progressions on a circle. For coprime integers $p\ge1$ and $r$, and $1\le h\le p$, let $g_1,\dots,g_h$ be the circular spacings of the points $0,r,\dots,(h-1)r$ in $\Z/p\Z$, and put
\[
\phi_{p,r}(h)=\sum_{i=1}^h\binom{g_i}2,\qquad\Phi(p,r)=\sum_{h=1}^p\phi_{p,r}(h),
\]
with $\phi_{p,r}(h)=0$ for $h>p$. Then $\phi_{p,r}(1)=\binom p2$, $\phi_{p,r}(p)=0$, and $\phi_{p,r}$ depends only on $r\bmod p$. We write $\M(p,q)=\M(\langle p,q\rangle)$ for coprime $p,q\ge1$.

\paragraph{Sumsets and undershoot.} Throughout, $A_0=\{0=a_0<a_1<\dots<a_k=d\}$ is a set of integers with $\gcd(A_0)=1$; dividing by the gcd reduces the general case to this one. For $h\ge0$, $hA_0$ is the set of sums of $h$ elements of $A_0$, that is, of at most $h$ elements of $A=A_0\setminus\{0\}$. For an integer $T\ge0$ the \emph{undershoot} is
\[
D_h(T)=T-\max\{n\in hA_0: n\le T\},
\]
and $\Tot(h)=\sum_{T=0}^{hd}D_h(T)$ is the total undershoot. We write $S_1=\langle a_1,\dots,a_k\rangle$ and $S_2=\langle d-a_0,\dots,d-a_{k-1}\rangle$ for the two numerical semigroups of $A_0$, with conductors $c_1$ and $c_2$.

\paragraph{The homogenized semigroup.} Let $\mathcal S\subset\N^2$ be the semigroup generated by the points $(a,d-a)$, $a\in A_0$. Then $n\in hA_0$ if and only if $(n,hd-n)\in\mathcal S$; the semigroup ring $\Bbbk[\mathcal S]$ over a field $\Bbbk$ is the homogeneous coordinate ring of the projective monomial curve $C_A$. Let $G$ be the group generated by $\mathcal S$.

\begin{theorem}[Cavaliere and Niesi {\cite[Lemma~4.3, Thm.~4.6]{CavaliereNiesi1983}}]\label{thm:cm}
The ring $\Bbbk[\mathcal S]$ is Cohen--Macaulay if and only if $\mathcal S=G\cap(S_1\times S_2)$. Equivalently, for every $h\ge1$,
\[
hA_0=I_h:=\{n\in[0,hd]: n\in S_1,\ hd-n\in S_2\}.
\]
\end{theorem}

We call $A_0$ \emph{CM} when these conditions hold, and $I_h$ the \emph{intersection form} at level $h$; in general $hA_0\subseteq I_h$. The same statement is \cite[Prop.~2.6]{GimenezGonzalezSanchez2023}; for general simplicial affine semigroups see \cite[Thm.~(2.6)]{GotoSuzukiWatanabe1976} and \cite[Cor.~4.4]{TrungHoa1986}.

\paragraph{Sumsets regularity.} Following \cite[Def.~1.3]{GimenezGonzalezSanchez2023}, the sumsets regularity $\sigma(A)$, a Hilbert-function index of $C_A$ and not its Castelnuovo--Mumford regularity, is the least $h_0\ge0$ such that for every $h\ge h_0$ the interval $[c_1,hd-c_2]$ is nonempty and
\[
hA_0=\bigl(S_1\cap[0,c_1-2]\bigr)\ \sqcup\ [c_1,hd-c_2]\ \sqcup\ \bigl(hd-(S_2\cap[0,c_2-2])\bigr).
\]
The requirement that the central interval be nonempty, that is $hd\ge c_1+c_2$, matters: without it the decomposition can hold at a smaller level (for $\{0,3,4,6\}$ it holds already at $h=1$ with an empty central interval). By \cite[Thm.~3.1, Prop.~3.2]{GimenezGonzalezSanchez2023}, $\sigma(A)$ is the maximum of the Hilbert regularity index of $\Bbbk[\mathcal S]$ and $\lceil(c_1+c_2)/d\rceil$, and equals $\lceil(c_1+c_2)/d\rceil$ when $A_0$ is CM. Gimenez and Gonz\'alez-S\'anchez assume that $A_0$ has at least four elements; for three elements $C_A$ is a plane curve of degree $d$, hence arithmetically Cohen--Macaulay with $\reg\Bbbk[\mathcal S]=d-1$, which is also the value of \eqref{eq:reg-cm} below, and the accompanying code checks the other facts cited from \cite{GimenezGonzalezSanchez2023} on the three-element sets of Table~\ref{tab:checks}.

We use four elementary facts throughout.

\begin{lemma}[Elementary facts]\label{lem:elementary}
\begin{enumerate}
\item[(i)] If $T\ge hd$ then $D_h(T)=T-hd$.
\item[(ii)] For parts with greatest common divisor $g$, the undershoot of a target $T$ is $(T\bmod g)+g\,D_h(\lfloor T/g\rfloor)$ for the reduced set.
\item[(iii)] Removing an element $x$ with $0<x<hd$ from $hA_0$ strictly increases $\Tot(h)$.
\item[(iv)] On each maximal run of $L$ consecutive non-members of $hA_0$ inside $[0,hd]$, the distances down to the preceding member and the distances up to the following member both add up to $L(L+1)/2$.
\end{enumerate}
\end{lemma}

\begin{proof}
(i) Every element of $hA_0$ is at most $hd$, and $hd\in hA_0$. (ii) The elements of $hA_0$ are the multiples $gn$ with $n$ in the $h$-fold sumset of the reduced set, so the largest one not above $T$ is $g$ times the largest one not above $\lfloor T/g\rfloor$. (iii) The undershoot at $x$ becomes positive and no undershoot decreases. (iv) Both sums are $1+2+\dots+L$.
\end{proof}

\section{Sumsets: residue intervals, counting and the spectrum}\label{sec:sumsets}

This section proves Theorem~\ref{thm:A}. It turns the intersection form of Theorem~\ref{thm:cm} into an exact count of the undershoot. The masses of $S_1$ and $S_2$ appear as the limits of the totals. A condition on two integer sequences, similar ordering, characterizes equality in the spectral lower bound and is sufficient, though not necessary, for the totals to be monotone: the Cohen--Macaulay set $\{0,2,3,5\}$ has $\Tot(h)=2$ for every $h\ge1$, while its only modulus, $5$, gives the incomparable pairs $(1,1)$ and $(2,0)$.

\subsection{Residue intervals and exact counting}\label{subsec:counting}

A \emph{modulus} for $A_0$ is a positive integer $m$ with $m\mid d$ and $m\in S_1\cap S_2$; $m=d$ is always one. Fix a modulus $m$, let $(w^{(1)}_r)$ and $(w^{(2)}_r)$ be the Ap\'ery sets of $S_1$ and $S_2$ with respect to $m$, and for $0\le s<m$ put
\[
u_s=\frac{s+w^{(1)}_{-s\bmod m}}{m},\qquad v_s=\frac{w^{(2)}_s-s}{m}.
\]
Subscripts of Ap\'ery elements are read modulo $m$. These are integers, $u_0=v_0=0$, $u_s\ge1$ for $s\ge1$ (as $w^{(1)}_{-s}>0$ is congruent to $m-s$), and $v_s\ge0$. For $h\ge0$ put $N=hd/m$.

\begin{lemma}[Residue intervals]\label{lem:residue}
Let $h\ge0$, $m$ a modulus, and $n=km-s$ with $0\le s<m$ and $0\le n\le hd$. Then $n\in I_h$ if and only if $u_s\le k\le N-v_s$. If $A_0$ is CM, the same holds for $hA_0$.
\end{lemma}

\begin{proof}
An integer $x\ge0$ lies in a numerical semigroup containing $m$ if and only if it is at least the Ap\'ery element of its class. The class of $n$ is $-s$, so $n\in S_1$ if and only if $km-s\ge w^{(1)}_{-s}$, that is $k\ge u_s$; and $hd-n=(N-k)m+s\in S_2$ if and only if $(N-k)m+s\ge w^{(2)}_s$, that is $k\le N-v_s$. For CM sets $hA_0=I_h$ (Theorem~\ref{thm:cm}; both are $\{0\}$ for $h=0$).
\end{proof}

For $m=d$ and CM sets this is the coordinate description of the homogenized Ap\'ery set in \cite[Prop.~2.6(3)]{GimenezGonzalezSanchez2023}, read on the line of level $h$. Every multiple $km$, $0\le k\le N$, lies in $I_h$, so every maximal run of integers of $[0,hd]$ outside $I_h$ lies in one \emph{block} $((k-1)m,km)$, $1\le k\le N$. Call $s\in[1,m-1]$ a \emph{hole} of block $k$ if $km-s\notin I_h$. By Lemma~\ref{lem:residue} the holes of block $k$ form
\begin{equation}\label{eq:holes}
H_k=B_k\cup R_{N-k},\qquad B_k=\{s\in[1,m-1]:u_s>k\},\quad R_j=\{s\in[1,m-1]:v_s>j\}.
\end{equation}

\begin{proposition}[Exact counting]\label{prop:counting}
Let $A_0$ be CM, $h\ge1$ and $m$ a modulus. For $1\le k\le N$ and $1\le s\le m-1$:
\begin{enumerate}
\item[(a)] $D_h(km-s)=\#\{t: s\le t\le m-1,\ [s,t]\subseteq H_k\}$;
\item[(b)] $\Tot(h)=\sum_{k=1}^{N}f(H_k)$;
\item[(c)] for integers $0\le T_-\le T_+\le hd$, $\sum_{T=T_-}^{T_+}D_h(T)$ is the number of triples $(k,s,t)$ with $1\le k\le N$, $1\le s\le t\le m-1$, $T_-\le km-s\le T_+$ and $[s,t]\subseteq H_k$.
\end{enumerate}
\end{proposition}

\begin{proof}
(a) If $s\notin H_k$ both sides vanish. Otherwise let $s'$ be the least non-hole of block $k$ above $s$, or $s'=m$. The predecessor of $km-s$ in $hA_0$ is $km-s'$, so $D_h(km-s)=s'-s$, the number of $t\in[s,s'-1]$, which are the $t$ with $[s,t]\subseteq H_k$. (b), (c) Every integer of $[0,hd]$ not divisible by $m$ is $km-s$ for exactly one such pair $(k,s)$, and multiples of $m$ have undershoot $0$; sum (a).
\end{proof}

In (c) only the target $km-s$ is restricted to the window, not $t$: the run below a target may leave the window. Lattice targets are handled by summing (a), targets $T\ge hd$ by Lemma~\ref{lem:elementary}(i) and common divisors by Lemma~\ref{lem:elementary}(ii). The proof uses CM only through $hA_0=I_h$, so the proposition holds for any $A_0$ at a level with $hA_0=I_h$, and for $I_h$ itself at every level.

\subsection{The spectrum and monotone totals}\label{subsec:spectrum}

Put $\beta_k=f(B_k)$ for $k\ge1$ and $\gamma_j=f(R_j)$ for $j\ge0$. Both sequences are nonincreasing and eventually $0$. The \emph{spectrum} of $A_0$ at the modulus $m$ is the multiset of all $\beta_k$ and $\gamma_j$, padded by zeros and written $\lambda_1\ge\lambda_2\ge\cdots$. We call $(u,v)$ \emph{similarly ordered} if $u_s>u_t$ implies $v_s\ge v_t$, that is, if the pairs $(u_s,v_s)$ form a chain in the product order; this is the hypothesis of Chebyshev's sum inequality and the equality case of the rearrangement inequality \cite[Thms.~43, 368 and 369]{HardyLittlewoodPolya1952}. Put $U=\max_su_s$ and $V=\max_sv_s$.

\begin{lemma}\label{lem:spectrum-mass}
$\sum_{k\ge1}\beta_k=\M(S_1)$ and $\sum_{j\ge0}\gamma_j=\M(S_2)$.
\end{lemma}

\begin{proof}
Multiples of $m$ lie in $S_1$ and $S_2$, so $\M(S_i)$ is the sum over blocks $k\ge1$ of $f$ of the gaps of $S_i$ in the block. As in Lemma~\ref{lem:residue}, the gaps of $S_1$ in block $k$ are the $km-s$ with $s\in B_k$, and a position $(k-1)m+s$, $1\le s<m$, lies in $S_2$ if and only if $k-1\ge v_s$, so the gaps of $S_2$ in block $k$ form a translate of $R_{k-1}$.
\end{proof}

\begin{lemma}[Merge]\label{lem:merge}
Let $(x_k)_{k\ge1}$, $(y_j)_{j\ge0}$ be nonincreasing sequences of nonnegative integers, eventually $0$, and $N\ge1$. Then $\sum_{k=1}^{N}\max(x_k,y_{N-k})$ is the sum of the $N$ largest entries of the multiset of all $x_k$ and $y_j$.
\end{lemma}

\begin{proof}
For $z\ge1$ let $p_z=\#\{k:x_k\ge z\}$ and $q_z=\#\{j:y_j\ge z\}$. Among $1\le k\le N$, $x_k\ge z$ for the first $\min(N,p_z)$ indices and $y_{N-k}\ge z$ for the last $\min(N,q_z)$, so $\max(x_k,y_{N-k})\ge z$ for $\min(N,p_z+q_z)$ values of $k$. The multiset has $p_z+q_z$ entries at least $z$. Sum over $z$.
\end{proof}

\begin{theorem}[Spectrum]\label{thm:spectrum}
Let $A_0$ be CM and $m$ a modulus.
\begin{enumerate}
\item For every $h\ge1$, $\Tot(h)\ge\lambda_1+\dots+\lambda_N$.
\item If $(u,v)$ is similarly ordered, equality holds in (1) for every $h$. For $m=d$, equality holds for every $h\ge1$ if and only if $(u,v)$ is similarly ordered.
\item If $(u,v)$ is similarly ordered, then $\Tot(h)=\sum_{\ell=1}^{hd/m}\lambda_\ell$ for every $h\ge0$. Hence $\Tot$ is nondecreasing and discretely concave, $\Tot(h)\le\M(S_1)+\M(S_2)$, and equality holds exactly when $N\ge\max(U-1,0)+V$, which for $m\ge2$ reads $N\ge U+V-1$.
\end{enumerate}
\end{theorem}

\begin{proof}
(1) The functional $f$ is monotone under inclusion, and strictly monotone under proper inclusion, since a point of $Y\setminus X$ is an interval counted by $f(Y)$ only. By Proposition~\ref{prop:counting}(b), $\Tot(h)=\sum_{k=1}^{N}f(B_k\cup R_{N-k})\ge\sum_{k=1}^{N}\max(\beta_k,\gamma_{N-k})$, and Lemma~\ref{lem:merge} identifies the right side.

(2) If $B_k$ and $R_j$ are not nested, there are $s\in B_k\setminus R_j$ and $t\in R_j\setminus B_k$, so $u_s>k\ge u_t$ and $v_s\le j<v_t$; under similar ordering this is impossible, and every block gives equality. Conversely let $m=d$, so $N=h$, and let $u_s>u_t$, $v_s<v_t$. Then $t\ne0$, so $k=u_t\ge1$; put $j=v_s$ and $h=k+j$. Then $s\in B_k\setminus R_j$ and $t\in R_j\setminus B_k$, so $f(B_k\cup R_j)>\max(\beta_k,\gamma_j)$ and (1) is strict at $h$.

(3) By (2) and Lemma~\ref{lem:merge}, $\Tot(h)$ is the sum of the $N=(d/m)h$ largest entries (also at $h=0$). Raising $h$ adds the next $d/m$ entries, which are nonnegative and entrywise at most those added before; with Lemma~\ref{lem:spectrum-mass} this gives the first claims. For $m\ge2$, $B_k\neq\varnothing$ exactly for $1\le k<U$ and $R_j\ne\varnothing$ exactly for $0\le j<V$, so exactly $U-1+V$ entries are positive. For $m=1$ there are no positions and every entry is $0$.
\end{proof}

In words, the bound takes in each block the larger of its two hole sets, one coming from $S_1$ and one from $S_2$, and it is exact in that block when one of the two contains the other. Similar ordering makes this happen in every block, and at $m=d$ nothing else does. For a general modulus the levels at which (1) is strict are determined by the inversions of $(u,v)$.

\begin{proposition}[Strict levels]\label{prop:strict-levels}
Let $A_0$ be CM, $m$ a modulus and $h\ge1$. Then $\Tot(h)>\lambda_1+\dots+\lambda_N$ if and only if
\[
N\in\bigcup_{u_s>u_t,\ v_s<v_t}[\,u_t+v_s,\ u_s+v_t-2\,].
\]
Thus equality holds at every level if and only if no positive multiple of $d/m$ lies in this finite union; for $m=d$ this is similar ordering.
\end{proposition}

\begin{proof}
By the proof of Theorem~\ref{thm:spectrum}, $\Tot(h)$ is strict exactly when some block $k$ has $B_k$, $R_{N-k}$ not nested, that is when there are $s,t$ with $u_t\le k\le u_s-1$ and $v_s\le N-k\le v_t-1$; for a fixed inversion $(s,t)$ such $k$ exists if and only if $N$ lies in the displayed interval (and then $1\le u_t\le k\le N$).
\end{proof}

\begin{example}\label{ex:spectrum}
(a) For $A_0=\{0,2,5,10\}$ and $m=10$ the positive spectrum is $(14,6,2,2)$; at $h=3$ the bound in Theorem~\ref{thm:spectrum}(1) is $22$, while $\Tot(3)=23$. At $m=5$ the pairs are similarly ordered (Lemma~\ref{lem:family-residue} in Section~\ref{sec:family}), so both notions depend on the modulus. (b) The CM set $\{0,2,8,11\}$ has $\M(S_1)+\M(S_2)=18$ and $\Tot(h)=19,19,18,18,\dots$: the totals are not monotone and exceed the mass. There are five inversions $(s,t)$, namely $(6,1)$, $(6,7)$, $(6,5)$, $(4,5)$ and $(10,5)$, with intervals $[1,1]$, $[1,1]$, $[1,2]$, $[2,2]$ and $[2,2]$, whose union is $[1,2]$, so by Proposition~\ref{prop:strict-levels} the inequality of Theorem~\ref{thm:spectrum}(1) is strict exactly at $h=1$ and $h=2$.
\end{example}

\subsection{The threshold}\label{subsec:threshold}

For a CM set and $m=d$ the homogenized Ap\'ery set consists of the points $(u_sd-s,v_sd+s)$, of level $u_s+v_s$ \cite[Prop.~2.6(3)]{GimenezGonzalezSanchez2023}, and the regularity of $\Bbbk[\mathcal S]$ is the largest level \cite[Thm.~3.7]{GimenezGonzalezSanchez2023}:
\begin{equation}\label{eq:reg-cm}
\reg \Bbbk[\mathcal S]=\max_{0\le s<d}(u_s+v_s).
\end{equation}
Under similar ordering one index attains both $U$ and $V$: if $u_s=U$ and $v_t=V$, then either $u_t=U$, or $u_s>u_t$ and $v_s\ge v_t$. So the maximum is $U+V$. (Similar ordering concerns the pairs, not the Ap\'ery points: for $\{0,1,2,5\}$ the indices $s=1,2$ both give $(1,1)$, while their points $(4,6)$ and $(3,7)$ are incomparable.)

\begin{proposition}[Threshold]\label{prop:threshold}
Let $A_0$ be CM and $m\ge2$ a modulus at which $(u,v)$ is similarly ordered. Then $U+V-1=\lceil(c_1+c_2)/m\rceil$, $\Tot(h)=\M(S_1)+\M(S_2)$ if and only if $h\ge\sigma(A)$, and $\reg \Bbbk[\mathcal S]=\sigma(A)+1$. If $m=d$, then
\[
U+V-1=\Bigl\lceil\frac{c_1+c_2}{d}\Bigr\rceil=\sigma(A).
\]
\end{proposition}

\begin{proof}
The Ap\'ery elements with respect to $m$ are $u_sm-s$ for $S_1$ and $v_sm+s$ for $S_2$. Let $s_1$ be the least $s$ with $u_s=U$ and $s_2$ the largest $s$ with $v_s=V$. Since $0\le s<m$, the largest Ap\'ery elements are $Um-s_1$ and $Vm+s_2$, so $c_1+c_2=(U+V-2)m+s_2-s_1+2$ (this also covers $S_i=\N$). As $m\ge2$, $U\ge1$ and $s_1\ge1$; an index attaining both maxima lies in $[s_1,s_2]$; so $2\le s_2-s_1+2\le m$, that is $(U+V-2)m<c_1+c_2\le(U+V-1)m$, and the ceiling is $U+V-1$. As $N=hd/m$ is an integer, $hd\ge c_1+c_2$ if and only if $N\ge U+V-1$. For CM sets $\sigma(A)=\lceil(c_1+c_2)/d\rceil$ \cite[Prop.~3.2]{GimenezGonzalezSanchez2023}, so $h\ge\sigma(A)$ if and only if $hd\ge c_1+c_2$, and Theorem~\ref{thm:spectrum}(3) gives the claim on $\Tot$. For $m=d$ we have $N=h$.

\emph{The regularity.} Let $s$ be an index with $u_s=U$ and $v_s=V$; then $s\ge1$, since $U\ge1$. Put $q=d/m$. Since $m\in S_1\cap S_2$, the Ap\'ery set of $S_i$ with respect to $d$ consists of the elements $w+jm$ with $w\in\Ap(S_i;m)$ and $0\le j<q$: each lies in $S_i$, and an element of $S_i$ congruent to $w+jm$ modulo $d$ is congruent to $w$ modulo $m$, hence equals $w+j'm$ with $j'\ge0$ and $j'\equiv j\pmod q$, so it is at least $w+jm$. Choose $0\le j_1,j_2<q$ with $q-1\le j_1+j_2\le2q-2$ such that $q$ divides $U+V+j_1+j_2$; the range holds $q$ consecutive integers, so exactly one value of $j_1+j_2$ qualifies. Then $x_1=(U+j_1)m-s\in\Ap(S_1;d)$, $x_2=(V+j_2)m+s\in\Ap(S_2;d)$ and $d$ divides $x_1+x_2$, so $x_1$ and $x_2$ are the Ap\'ery elements with respect to $d$ of the classes $-s'$ and $s'$, where $s'=x_2\bmod d$, and the pair of residue $s'$ at the modulus $d$ has sum $(x_1+x_2)/d$. As $U+V+j_1+j_2$ is the least multiple of $q$ that is at least $U+V+q-1$, and $\lceil\lceil x/m\rceil/q\rceil=\lceil x/d\rceil$ for $x\ge0$, \eqref{eq:reg-cm} gives
\[
\reg \Bbbk[\mathcal S]\ \ge\ \frac{U+V+j_1+j_2}{q}=\Bigl\lceil\frac{U+V-1}{q}\Bigr\rceil+1=\Bigl\lceil\frac{c_1+c_2}{d}\Bigr\rceil+1=\sigma(A)+1 .
\]
The reverse inequality is \cite[Thm.~3.16]{GimenezGonzalezSanchez2023}.
\end{proof}

So under similar ordering the total reaches the mass at the level from which $hA_0$ keeps its final shape.

For CM sets with $m=1$ both semigroups are $\N$ and every total vanishes. For $d=1$, $A_0=\{0,1\}$, $\sigma(A)=0$, $\Bbbk[\mathcal S]$ is a polynomial ring of regularity $0$, and both masses and all totals vanish. The nonempty central interval in the definition of $\sigma(A)$ is needed: $\{0,3,4,6\}$ is CM and similarly ordered with $(c_1,c_2)=(6,2)$, and its decomposition with an empty central interval holds already at $h=1<\lceil8/6\rceil$. Similarly ordered CM sets with $d\ge2$ attain the bound $\reg \Bbbk[\mathcal S]\le\sigma(A)+1$ of \cite[Thm.~3.16]{GimenezGonzalezSanchez2023}. The proof uses similar ordering only through an index attaining both $U$ and $V$, so every CM set with such an index at some modulus $m\ge2$ has $U+V-1=\lceil(c_1+c_2)/m\rceil$ and $\reg \Bbbk[\mathcal S]=\sigma(A)+1$. The hypothesis of \cite[Prop.~3.18]{GimenezGonzalezSanchez2023} amounts for CM sets to $s_1=s_2$, which gives such an index. It does not imply similar ordering: $\{0,2,5,10\}$ has $s_1=s_2=7$ at $m=10$ and is not similarly ordered there (Example~\ref{ex:spectrum}).

\section{Similarly ordered four-element sets}\label{sec:similar}\label{subsec:similar}

This section proves Theorem~\ref{thm:B}. For four elements the residue pairs of Section~\ref{sec:sumsets} at $m=d$ come from a two-dimensional diagram of $d$ lattice points, and similar ordering becomes a condition on its layers that can be checked from the two Ap\'ery sets.

Let $A_0=\{0,a_1,a_2,d\}$, $0<a_1<a_2<d$. Let $w^{(1)}_r$ and $w^{(2)}_r$ be the Ap\'ery elements of $S_1$ and $S_2$ modulo $d$, as in Section~\ref{subsec:counting} with $m=d$ (subscripts modulo $d$), and for $0\le r<d$ put
\begin{equation}\label{eq:diagram}
t_r=\frac{w^{(1)}_r+w^{(2)}_{-r}}{d},\qquad y_r=\frac{w^{(1)}_r-a_1t_r}{a_2-a_1},\qquad x_r=t_r-y_r .
\end{equation}

\begin{lemma}[Common diagram]\label{lem:diagram}
A four-element set is CM if and only if all $x_r,y_r$ are nonnegative integers. Then $\mathcal D=\{(x_r,y_r):0\le r<d\}$ consists of $d$ points, it is an order ideal of $\N^2$, the point $(x,y)$ of residue $r$ satisfies $w^{(1)}_r=a_1x+a_2y$ and $w^{(2)}_{-r}=(d-a_1)x+(d-a_2)y$, and with $s=-r\bmod d$,
\begin{equation}\label{eq:uv-diagram}
u_s=\Bigl\lceil\frac{a_1x+a_2y}{d}\Bigr\rceil,\qquad v_s=x+y-u_s .
\end{equation}
\end{lemma}

\begin{proof}
The point $(w^{(1)}_r,w^{(2)}_{-r})$ lies in $G\cap(S_1\times S_2)$ and is the coordinatewise least point there with first coordinate $\equiv r$. If $A_0$ is CM it lies in $\mathcal S$ (Theorem~\ref{thm:cm}); a representation cannot use $(d,0)$ or $(0,d)$, by minimality of $w^{(1)}_r$ and $w^{(2)}_{-r}$, so it is $x(a_1,d-a_1)+y(a_2,d-a_2)$, and solving gives $(x,y)=(x_r,y_r)$. Conversely, if the $x_r,y_r$ are nonnegative integers, every point of $G\cap(S_1\times S_2)$ is some $(w^{(1)}_r,w^{(2)}_{-r})\in\mathcal S$ plus multiples of $(d,0)$ and $(0,d)$, so $A_0$ is CM. If $(x',y')\le(x,y)\in\mathcal D$ and an element of $S_1$ smaller than $a_1x'+a_2y'$ were in its class, adding $a_1(x-x')+a_2(y-y')$ would contradict the minimality of $a_1x+a_2y$; the same holds for $S_2$, so $(x',y')\in\mathcal D$. Finally $u_s=(d-r+w^{(1)}_r)/d=\lceil w^{(1)}_r/d\rceil$ for $r>0$, and $u_s+v_s=t_r$.
\end{proof}

For $(x,y)\in\N^2$ call $a_1x+a_2y$ and $(d-a_1)x+(d-a_2)y$ its two values and $a_1x+a_2y\bmod d$ its residue. For a nonempty layer $\mathcal D_t=\{(x,y)\in\mathcal D:x+y=t\}$ let $F^-_t$ and $F^+_t$ be the least and largest values of $a_1x+a_2y=a_1t+(a_2-a_1)y$ on it, attained at the least and largest $y$.

\begin{theorem}[Similar ordering]\label{thm:similar}
For a CM set $\{0,a_1,a_2,d\}$ the following are equivalent:
(1) $(u,v)$ is similarly ordered at $m=d$;
(2) $\lceil(a_1x+a_2y)/d\rceil$ is constant on each layer $\mathcal D_t$;
(3) $\lceil F^-_t/d\rceil=\lceil F^+_t/d\rceil$ for every nonempty layer, that is, each layer is mapped into one interval $((q-1)d,qd]$ with $q\in\Z$.
\end{theorem}

\begin{proof}
By \eqref{eq:uv-diagram} two points of $\mathcal D_t$ give pairs $(u,t-u)$, $(u',t-u')$, incomparable if $u\ne u'$; so (1) implies (2). Assume (2), with value $q_t$ on $\mathcal D_t$. The nonempty layers are $0,\dots,t^*$, and every point of layer $t+1$ lies above a point of layer $t$ whose value $a_1x+a_2y$ is smaller by $a_1$ or $a_2\in(0,d)$, so $q_{t+1}-q_t\in\{0,1\}$. Hence $q_t$ and $t-q_t$ are nondecreasing and the pairs $(q_t,t-q_t)$ form a chain. (2) and (3) are equivalent because the value increases with $y$ on a layer.
\end{proof}

\begin{proposition}[Shape]\label{prop:shape}
For a CM set $\{0,a_1,a_2,d\}$, $\mathcal D$ is a rectangle $[0,P-1]\times[0,Q-1]$ or an L-shape $([0,W-1]\times[0,Q-1])\cup([W,P-1]\times[0,Z-1])$ with $0<W<P$, $0<Z<Q$, and then $a_1W+a_2Z\equiv0\pmod d$. If $P-W\ge2$ and $Q-Z\ge2$, $(u,v)$ is not similarly ordered.
\end{proposition}

\begin{proof}
The minimal points of $\N^2\setminus\mathcal D$ are $(P,0)$, $(0,Q)$ and possibly points $(p,q)$ with $p,q>0$. For such a point let $(x,y)\in\mathcal D$ be the point of its residue; its two values are at most those of $(p,q)$, one strictly. If $x>0$, then $(p-1,q)\in\mathcal D$ while $(x-1,y)$ has the same residue and values at most those of $(p-1,q)$, one strictly, contradicting Lemma~\ref{lem:diagram}; similarly if $y>0$. So $(x,y)=(0,0)$ and $a_1p+a_2q\equiv0$. Two such points $(p_1,q_1)$, $(p_2,q_2)$ with $p_1<p_2$ would put the distinct points $(p_2-p_1,0)$ and $(0,q_1-q_2)$ of the same residue in $\mathcal D$ (as $p_2-p_1<P$ and $q_1-q_2<Q$). If $P-W,Q-Z\ge2$, the points $(W\mp1,Z\pm1)$ lie on one layer with values $a_1W+a_2Z\pm(a_2-a_1)$ around a positive multiple of $d$, and Theorem~\ref{thm:similar} applies.
\end{proof}

\begin{figure}[ht]
\centering
\begin{tikzpicture}[x=0.88cm,y=0.88cm,every node/.style={font=\scriptsize}]
\begin{scope}
\begin{scope}\clip (-0.45,-0.45) rectangle (6.6,1.55);
\foreach \t in {1,...,7} {\draw[gray,dashed] (\t+0.45,-0.45) -- (-0.45,\t+0.45);}
\end{scope}
\foreach \x in {0,...,6} {\foreach \y in {0,1} {
  \pgfmathtruncatemacro{\v}{5*\x+7*\y}
  \fill (\x,\y) circle (1.6pt); \node[above right=-1pt] at (\x,\y) {\v};}}
\node at (3,-0.9) {$\{0,5,7,14\}$: rectangle, similarly ordered};
\end{scope}
\begin{scope}[xshift=7.6cm]
\begin{scope}\clip (-0.45,-0.45) rectangle (3.6,2.55);
\foreach \t in {1,4} {\draw[gray,dashed] (\t+0.45,-0.45) -- (-0.45,\t+0.45);}
\foreach \t in {2,3} {\draw[black,thick] (\t+0.45,-0.45) -- (-0.45,\t+0.45);}
\end{scope}
\foreach \x/\y in {0/0,1/0,2/0,3/0,0/1,1/1,2/1,3/1,0/2,1/2,2/2} {
  \pgfmathtruncatemacro{\v}{2*\x+8*\y}
  \fill (\x,\y) circle (1.6pt); \node[above right=-1pt] at (\x,\y) {\v};}
\node at (1.5,-0.9) {$\{0,2,8,11\}$: L-shape, not similarly ordered};
\end{scope}
\end{tikzpicture}
\caption{Ap\'ery diagrams, each point $(x,y)$ labeled with $a_1x+a_2y$; the lines are the layers $x+y=t$. For $\{0,5,7,14\}$ the values on every layer lie in one interval $(14(q-1),14q]$. For $\{0,2,8,11\}$ the layers $t=2$ and $t=3$ (solid) carry the values $4,10,16$ and $6,12,18$, which cross $11$. Drawn with Ti\emph{k}Z code written by Claude Opus (Anthropic) from the computed Ap\'ery sets.}\label{fig:diagrams}
\end{figure}
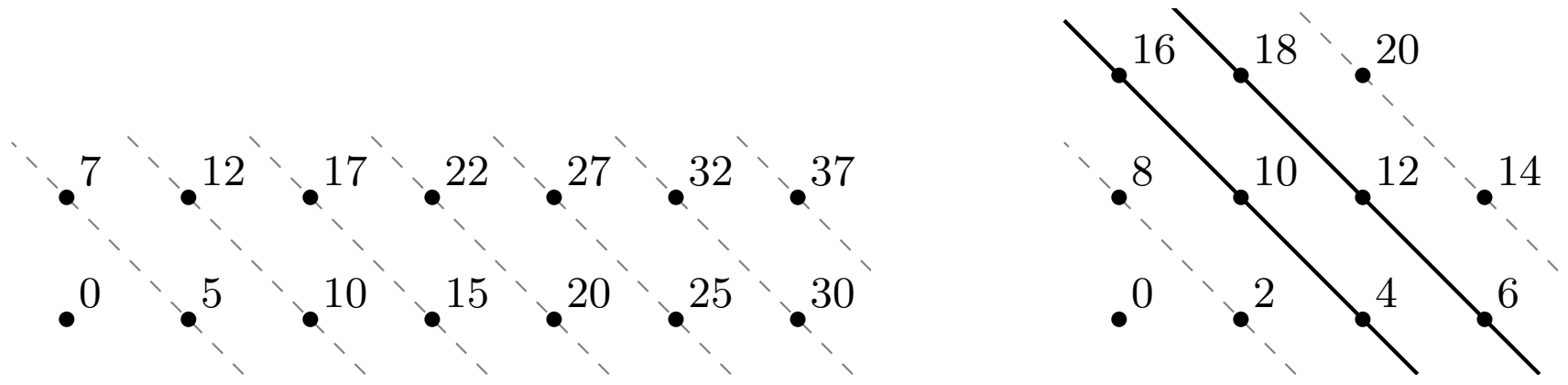

That $\mathcal D$ is a rectangle or an L-shape is the classical shape of the Ap\'ery set of the three-generated semigroup $S_1=\langle a_1,a_2,d\rangle$ with respect to $d$, drawn in the coordinates of $a_1$ and $a_2$ \cite{AguiloMarijuan2014}; what is used here is the congruence of the notch and its consequence for similar ordering. Figure~\ref{fig:diagrams} shows the test of Theorem~\ref{thm:similar} on two sets of the introduction.

\begin{proposition}[Families]\label{prop:families}
\begin{enumerate}
\item[(a)] Let $1<b<d$ and $d=qb+e$, $0\le e<b$. Then $\{0,1,b,d\}$ is CM if and only if $e=0$ or $q+e\ge b$, and then it is similarly ordered.
\item[(b)] $A_0$ is CM and similarly ordered if and only if $d-A_0$ is. In particular every CM set $\{0,a,d-1,d\}$ is similarly ordered.
\item[(c)] A CM set $\{0,a,d-a,d\}$, $0<a<d-a$, is similarly ordered if and only if $d=3$.
\item[(d)] Every set $\{0,a,a+1,d\}$, $a+1<d$, is CM, and it is similarly ordered if and only if $\mathcal D$ is a rectangle or an L-shape with $P-W=1$ or $Q-Z=1$.
\end{enumerate}
\end{proposition}

\begin{proof}
For $n\in\langle a_1,a_2\rangle$ let $\operatorname{len}(n)$ be the least number of parts $a_1,a_2$ with sum $n$. By the proof of Lemma~\ref{lem:diagram}, $A_0$ is CM if and only if $w^{(2)}_{-r}=\operatorname{len}(w^{(1)}_r)d-w^{(1)}_r$ for all $r$. The elements of $S_2$ of class $-r$ avoiding the part $d$ are $\ell'd-(w^{(1)}_r+kd)$ for representations of $w^{(1)}_r+kd$ by $\ell'$ interior parts, with $k\ge0$ by minimality of $w^{(1)}_r$. Hence $A_0$ is CM if and only if
\begin{equation}\label{eq:lambda-test}
\operatorname{len}(w^{(1)}_r+kd)-\operatorname{len}(w^{(1)}_r)\ge k\qquad(0\le r<d,\ k\ge0,\ w^{(1)}_r+kd\in\langle a_1,a_2\rangle).
\end{equation}
(a) Here $\Ap(S_1;d)=[0,d-1]$, and the shortest representation of $n$ uses $\lfloor n/b\rfloor$ parts $b$ and $n\bmod b$ parts $1$, so $\operatorname{len}(n)=\lfloor n/b\rfloor+(n\bmod b)$. Adding $d=qb+e$ raises $\operatorname{len}$ by $q$ if $e=0$. Otherwise it raises it by $q+e$, or by $q+e-b+1$ when the remainder passes $b$. If $q+e\ge b$, every increment is at least $1$ and \eqref{eq:lambda-test} follows. If $0<e$ and $q+e<b$, the element $n=b-e$ of $\Ap(S_1;d)$ violates it, since $n+d=(q+1)b$ and its length rises by $q+e-b+1\le0$. Since $w^{(1)}_r=r$, we have $u_s=1$ for every $s\ge1$, and the pairs $(1,v_s)$ form a chain together with $(0,0)$.
(b) Reflection exchanges the two coordinates of $\mathcal S$, so it preserves the CM property; it exchanges $S_1$ and $S_2$, and with $s'=d-s$ the reflected pairs are $(v_s+1,u_s-1)$ for $s\ne0$. Swapping the coordinates of two pairs and shifting them by the same amounts keeps them comparable or incomparable, so similar ordering is preserved. The second claim follows from (a).
(c) Here reflection fixes the set, so with each pair $(u,v)$, $s\ne0$, the pair $(v+1,u-1)$ occurs on the same layer; comparability forces $u=v+1$, so the layer $u+v$ is odd. For $d>3$ the order ideal $\mathcal D$ has more than three points, while layers $0$ and $1$ hold at most three, hence a point on layer $2$, which is impossible. For $d=3$ the pairs are $(0,0),(1,0),(1,0)$.
(d) Let $b=a+1$. The representation of $n\in\langle a,b\rangle$ with fewest parts has $x(n)=(-n)\bmod b\le a$ parts $a$, so $\operatorname{len}(n)=(n+x(n))/b$, and for $k\ge1$, $\operatorname{len}(n+kd)-\operatorname{len}(n)\ge(kd-a)/b>k-1$ as $d>b$; so \eqref{eq:lambda-test} holds. If both notch dimensions are at least $2$, use Proposition~\ref{prop:shape}. Otherwise each layer is a set of consecutive lattice points (on an L-shape the points of layer $t$ with $x<W$ and those with $x\ge W$ are separated only if $(W,t-W)$ lies in the notch while both groups are nonempty, that is $W+Z\le t\le\min(W+Q-2,P+Z-2)$, a range that is empty when $P-W=1$ or $Q-Z=1$), with consecutive values $at+y$; none of these is a multiple of $d$ for $t>0$, since it would share the residue of $(0,0)$, so the layer lies in one interval $((q-1)d,qd)$.
\end{proof}

The condition in (a) is also the condition under which the cardinality formula of Chen and Yang for $|h\{0,1,b,d\}|$ holds \cite[Cor.~1]{ChenYang2026}, obtained there by the same count of least numbers of parts; the similar-ordering conclusion is not in their statement.

The family $\{0,a,b,2b\}$ of Section~\ref{sec:family} is similarly ordered at $m=b$ (Lemma~\ref{lem:family-residue}) but not always at $m=2b$: $\{0,2,5,10\}$ has $\mathcal D=[0,4]\times[0,1]$, and its layer $4$ has values $8$ and $11$. Theorem~\ref{thm:similar} is a finite arithmetic test from the two Ap\'ery sets; a description of the similarly ordered triples $(a_1,a_2,d)$ by a short list of congruence or continued-fraction conditions is not known (Section~\ref{sec:open}).

\section{The family \texorpdfstring{$\{0,a,b,2b\}$}{\{0,a,b,2b\}}}\label{sec:family}

This section proves Theorem~\ref{thm:C}. At the modulus $b$ the pairs of the family are similarly ordered, its missing positions are complements of arithmetic progressions on $\Z/b\Z$, and the spectrum becomes a sum over circular gaps, which the three-distance theorem evaluates.

Throughout this section $0<a<b$ are coprime and $A_0=\{0,a,b,2b\}$, so $d=2b$, $S_1=\langle a,b\rangle$ and $S_2=\langle b,2b-a\rangle$. We put $\M_4(a,b)=\M(a,b)+\M(b,2b-a)$ and $\phi(j)=\phi_{b,a}(j)$ (Section~\ref{sec:prelim}), so $\phi(1)=\binom b2$ and $\phi(j)=0$ for $j\ge b$.

\subsection{Structure and spectrum}\label{subsec:family-structure}

\begin{lemma}\label{lem:family-ci}
Over every field, the ideal of $C_A$ is generated by $x_2^2-x_0x_3$ and $x_1^b-x_2^{\,j}x_0^{\,z}x_3^{\,w}$, where $x_0,\dots,x_3$ correspond to the parts $0,a,b,2b$, $j=a\bmod2$, $w=(a-j)/2$ and $z=(2b-a-j)/2$. So $C_A$ is a complete intersection of degrees $2$ and $b$ and $A_0$ is CM; $|hA_0|=(h+1)^2$ for $h\le b-1$ and $|hA_0|=2bh+1-(b-1)^2$ for $h\ge b-1$; $c_1=(a-1)(b-1)$, $c_2=(b-1)(2b-a-1)$ and $\sigma(A)=b-1$.
\end{lemma}

\begin{proof}
Both binomials are homogeneous (of degrees $2$ and $b=j+w+z$) and lie in the ideal $I$ of $C_A$. Modulo them every monomial is congruent to a monomial $x_0^rx_1^ex_2^\varepsilon x_3^t$ of the same degree with $e<b$, $\varepsilon\le1$: replacing $x_1^b$ lowers the exponent of $x_1$, replacing $x_2^2$ keeps it and lowers that of $x_2$. Such a monomial maps to the point with first coordinate $ae+b\varepsilon+2bt$; for two of them of equal degree and equal image, $ae\equiv ae'\pmod b$ gives $e=e'$, then parity gives $\varepsilon=\varepsilon'$ and $t=t'$, and the degree gives $r=r'$. So in each degree these monomials have linearly independent images, every homogeneous element of $I$ reduces to $0$, and $I$ is generated by the two binomials; a prime of height two with two generators is a complete intersection (compare \cite{BermejoGarciaMarco2014}). Hence $|hA_0|$ is the number of these monomials of degree $h$, namely $\sum_{e=0}^{\min(h,b-1)}(2(h-e)+1)$. The conductors are those of $\langle a,b\rangle$ and $\langle b,2b-a\rangle$, $c_1+c_2=2(b-1)^2$, and $\sigma(A)=\lceil(b-1)^2/b\rceil=b-1$ (Section~\ref{sec:prelim}).
\end{proof}

The integer $b$ is a modulus (Section~\ref{subsec:counting}), with $N=2h$.

\begin{lemma}[Residue intervals modulo $b$]\label{lem:family-residue}
For $0\le s<b$ let $i_s\in[0,b-1]$ satisfy $i_sa\equiv-s\pmod b$. At $m=b$, $u_s=(s+i_sa)/b=\lceil i_sa/b\rceil$ and $v_s=2i_s-u_s$, so $kb-s\in hA_0$ if and only if $u_s\le k\le2h-v_s$. As functions of $i=i_s$, $u$ is nondecreasing with steps $0,1$ and $v$ is increasing with steps $1,2$. Hence $(u,v)$ is similarly ordered at $m=b$, $U=a$, $V=2b-a-2$, and $\Tot$ is nondecreasing and concave and equals $\M_4(a,b)$ exactly for $h\ge b-1$.
\end{lemma}

\begin{proof}
The Ap\'ery elements modulo $b$ are $ia$ for $S_1$ and $i(2b-a)\equiv-ia$ for $S_2$, so $w^{(1)}_{-s}=i_sa$ and $w^{(2)}_s=i_s(2b-a)$, which gives $u_s$ and $v_s$; Lemma~\ref{lem:residue} applies by Lemma~\ref{lem:family-ci}. Since $0<a/b<1$, $\lceil(i+1)a/b\rceil-\lceil ia/b\rceil\in\{0,1\}$. If $u_s>u_t$ then $i_s>i_t$, because $u$ is nondecreasing in $i$, and then $v_s>v_t$, because $v$ is increasing. Both maxima occur at $i=b-1$. By Theorem~\ref{thm:spectrum}(3), the mass is reached exactly when $2h\ge U+V-1=2b-3$.
\end{proof}

The family need not be similarly ordered at $d=2b$ (Section~\ref{subsec:similar}); at the modulus $b$, Proposition~\ref{prop:threshold} gives the level $\sigma(A)=b-1$ of Lemma~\ref{lem:family-ci} again, and $\reg\Bbbk[\mathcal S]=b$, the regularity of a complete intersection of degrees $2$ and $b$.

For $1\le i\le b-1$ let $s(i)\in[1,b-1]$ satisfy $s(i)\equiv-ia\pmod b$, put $\mathcal P_j=\{s(i):j\le i\le b-1\}$, and let $\mathcal R_j$ be the multiset of the numbers $g-1\ge1$ over the gaps $g$ of $\{0,a,\dots,(j-1)a\}$ on $\Z/b\Z$.

\begin{theorem}[Run lengths]\label{thm:family-runs}
For $h\ge1$ the lengths of the maximal runs of integers of $[0,2hb]$ outside $hA_0$ form the multiset
\[
\mathcal R_1\sqcup\bigsqcup_{j=2}^{h}2\mathcal R_j\sqcup\mathcal R_{h+1}\quad(h<b-1),\qquad \mathcal R_1\sqcup\bigsqcup_{j=2}^{b-1}2\mathcal R_j\quad(h\ge b-1),
\]
where $2\mathcal R_j$ is two copies. In particular the distribution of $D_h$ on $[0,2hb]$ is unchanged when $a$ is replaced by $b-a$.
\end{theorem}

\begin{proof}
By Lemma~\ref{lem:family-residue}, position $s(i)$ is present in block $k$ exactly for $k\in E_i=[\alpha_i,\alpha_i+2h-2i]$, $\alpha_i=\lceil ia/b\rceil$, an interval inside $[1,2h]$ when nonempty. As $\alpha_{i+1}-\alpha_i\in\{0,1\}$, $E_{i+1}\subseteq E_i$, and $|E_i|=\max(2(h-i)+1,0)$. So the holes of block $k$ are $\mathcal P_j$ with $j=1+\max(\{0\}\cup\{i:k\in E_i\})$. Since $E_1=[1,2h-1]$, one block has holes $\mathcal P_1$. For $2\le j\le\min(h,b-1)$ the two blocks of $E_{j-1}\setminus E_j$ have holes $\mathcal P_j$. If $h<b-1$, the block of $E_h$ has holes $\mathcal P_{h+1}$, and if $h\ge b-1$, the blocks of $E_{b-1}$ have none. The complement of $\mathcal P_j$ in $[1,b-1]$, together with $0$, is $\{0,-a,\dots,-(j-1)a\}$, and a gap $g$ of this set is a run of $g-1$ holes. The map $x\mapsto-x$ carries the gaps of $\{0,a,\dots,(j-1)a\}$ to those of this set, so the runs of $\mathcal P_j$ have lengths $\mathcal R_j$. Replacing $a$ by $b-a$ reflects each point set and preserves each $\mathcal R_j$, and a run of length $L$ carries the undershoots $1,\dots,L$.
\end{proof}

\begin{theorem}[Spectrum of the family]\label{thm:family-spectrum}
At $m=b$ the spectrum of $A_0$ is $\binom b2$ once and $\phi(j)$ twice for $2\le j\le b-1$, and $\phi$ is nonincreasing. Hence
\begin{equation}\label{eq:family-tot}
\Tot(h)=\binom b2+2\sum_{j=2}^{h}\phi(j)+\phi(h+1)\qquad(1\le h\le b-1),
\end{equation}
and $\Tot(h)=\M_4(a,b)$ for $h\ge b-1$.
\end{theorem}

\begin{proof}
Write $u(i)=\lceil ia/b\rceil$ and $v(i)=2i-u(i)$ for the pairs of Lemma~\ref{lem:family-residue} as functions of $i$. Then $B_k=\mathcal P_j$ for the $u(j)-u(j-1)$ values $u(j-1)\le k<u(j)$, and $R_\ell=\mathcal P_j$ for the $v(j)-v(j-1)$ values $v(j-1)\le\ell<v(j)$. The two counts add up to $2$ for $j\ge2$, while $\mathcal P_1$ occurs only as $R_0$, since $k\ge1=u(1)$. By the proof of Theorem~\ref{thm:family-runs}, $f(\mathcal P_j)=\phi(j)$. Adding a point to a circular set splits a gap $x+y$ and lowers $\binom{x+y}2$ by $xy\ge0$, so $\phi$ is nonincreasing, and Theorem~\ref{thm:spectrum}(3) gives $\Tot(h)$ as the sum of the $2h$ largest entries.
\end{proof}

\begin{corollary}[Mass of the family]\label{cor:family-mass}
For coprime $0<a<b$, $\M_4(a,b)=2\Phi(b,a)-\binom b2$.
\end{corollary}

\begin{proof}
By Lemma~\ref{lem:spectrum-mass} the mass is the sum of the spectrum, $\binom b2+2\sum_{j=2}^{b-1}\phi(j)=2\Phi(b,a)-\binom b2$, as $\phi(1)=\binom b2$ and $\phi(b)=0$.
\end{proof}

\subsection{The whole range through continued fractions}\label{subsec:family-cf}

Let $a/b=[0;\theta_1,\dots,\theta_K]$ with $\theta_K\ge2$, with convergents $p_j/q_j$, $p_{-1}=1$, $q_{-1}=0$, $p_0=0$, $q_0=1$, so $q_K=b$, and put $\varepsilon_j=|q_ja-p_jb|$ for $-1\le j\le K$. Then $\varepsilon_{-1}=b$, $\varepsilon_0=a$, $\varepsilon_{K-1}=1$, $\varepsilon_K=0$ and $\varepsilon_{j-1}=\theta_{j+1}\varepsilon_j+\varepsilon_{j+1}$. The intervals $[\nu q_k+q_{k-1},(\nu +1)q_k+q_{k-1})$, $1\le \nu \le \theta_{k+1}$, tile $[q_k+q_{k-1},q_{k+1}+q_k)$, and these tile $[1,b+q_{K-1})$ for $0\le k<K$. So each $h\in[2,b]$ has a unique expression
\begin{equation}\label{eq:cf-index}
h-1=\nu q_k+q_{k-1}+r,\qquad 0\le k\le K-1,\ 1\le \nu \le \theta_{k+1},\ 0\le r<q_k ,
\end{equation}
and the block $(k,\nu )=(K-1,\theta_K)$, which would start at $h=b+1$, does not occur. The $h$ with a given $(k,\nu )$ form a \emph{block}.

The next theorem is the rational form of the classical three-distance theorem, in the notation above.

\begin{theorem}[Three distances]\label{thm:three-distance}
For $2\le h\le b$ as in \eqref{eq:cf-index}, the gaps of $\{0,a,\dots,(h-1)a\}$ on $\Z/b\Z$ are $h-q_k$ gaps $\varepsilon_k$, $r+1$ gaps $\varepsilon_{k-1}-\nu \varepsilon_k$ and $q_k-r-1$ gaps $\varepsilon_{k-1}-(\nu -1)\varepsilon_k$. Hence
\[
\phi(h)=(h-q_k)\binom{\varepsilon_k}2+(r+1)\binom{\varepsilon_{k-1}-\nu \varepsilon_k}2+(q_k-r-1)\binom{\varepsilon_{k-1}-(\nu -1)\varepsilon_k}2 .
\]
\end{theorem}

\begin{proof}
For irrational $\alpha$ with convergents $p_j/q_j$ and $\eta_j=|q_j\alpha-p_j|$, the three distance theorem in the form of \cite[Sect.~3]{AlessandriBerthe1998} (due to S\'os \cite{Sos1958}, Sur\'anyi \cite{Suranyi1958} and \'Swierczkowski \cite{Swierczkowski1958}) states that for $n=\nu q_k+q_{k-1}+r$ as above the points $0,\alpha,\dots,n\alpha$ cut $\mathbb R/\Z$ into $n+1-q_k$ arcs $\eta_k$, $r+1$ arcs $\eta_{k-1}-\nu \eta_k$ and $q_k-r-1$ arcs $\eta_{k-1}-(\nu -1)\eta_k$. Let $\alpha_N=[0;\theta_1,\dots,\theta_K,N,1,1,\dots]$: its convergents of index at most $K$ are those of $a/b$, and $\alpha_N\to a/b$. Take $n=h-1$; its expression involves only $k\le K-1$ and $\nu \le \theta_{k+1}$, so it is the same for $\alpha_N$. The points $ia/b\bmod1$, $0\le i<h\le b$, are distinct with circular spacing at least $1/b$. Once $|\alpha_N-a/b|<1/(3b^2)$, each $i\alpha_N$ lies within $1/(3b)$ of $ia/b$, so the circular order is the same and each arc converges to the corresponding rational gap divided by $b$, while $\eta_j(\alpha_N)\to\varepsilon_j/b$ for $j\le K-1$. Pass to the limit and multiply by $b$.
\end{proof}

Thus the formula for $\phi$ changes only from one block to the next. On a block, raising $h$ splits a gap $x+\varepsilon_k$, $x=\varepsilon_{k-1}-\nu \varepsilon_k$, into $x$ and $\varepsilon_k$, so $\phi$ is affine there with slope $-\delta_{k,\nu }$, $\delta_{k,\nu }=\varepsilon_k(\varepsilon_{k-1}-\nu \varepsilon_k)$. The quadratic pieces of $\Tot$ are shifted.

\begin{corollary}[Quadratic pieces]\label{cor:quadratic-pieces}
Put $\Tot(0)=0$, and let $h^*=\nu q_k+q_{k-1}+1$ be the first $h$ of block $(k,\nu )$. Then
\[
\Tot(h^*-2+j)=\Tot(h^*-2)+2\phi(h^*-1)\,j-\delta_{k,\nu }\,j^2\qquad(0\le j\le q_k),
\]
and the intervals $[h^*-2,h^*+q_k-2]$ tile $[0,b-1]$ with shared endpoints.
\end{corollary}

\begin{proof}
The formula of Theorem~\ref{thm:three-distance} at $r=-1$ equals $\phi(h^*-1)$. For $\nu >1$ it is the value of block $(k,\nu -1)$ at $r=q_k-1$. For $\nu =1$ and $k\ge1$ it is $q_{k-1}\binom{\varepsilon_k}2+q_k\binom{\varepsilon_{k-1}}2$, the value at the last $h=q_k+q_{k-1}$ of block $(k-1,\theta_k)$, because $\varepsilon_{k-2}-\theta_k\varepsilon_{k-1}=\varepsilon_k$. For $k=0$ and $\nu =1$ both are $\binom b2$. So $\phi(h^*-1+i)=\phi(h^*-1)-\delta_{k,\nu }i$ for $0\le i\le q_k$. By \eqref{eq:family-tot}, $\Tot(h)-\Tot(h-1)=\phi(h)+\phi(h+1)$ for $1\le h\le b-1$, which at $h=h^*-2+j$, $1\le j\le q_k$, is $2\phi(h^*-1)-\delta_{k,\nu }(2j-1)$, and summing over $j$ gives the formula. Consecutive blocks start at $h^*$ and $h^*+q_k$, also across rows, as $\theta_{k+1}q_k+q_{k-1}+q_k=q_{k+1}+q_k$. The first block has $h^*=2$ and the last ends at $h=b$.
\end{proof}

For $a=2$, $b=9$ ($[0;4,2]$) the last block of $\phi$ is $h=6,\dots,9$, but $\Tot(6),\dots,\Tot(9)=132,135,136,136$ has third difference $1$; the last quadratic piece of $\Tot$ is $4\le h\le8$, with values $120,127,132,135,136$.

\section{Non-Cohen--Macaulay sets}\label{sec:noncm}

This section proves Theorem~\ref{thm:D}. A set that is not CM has fewer sums than its intersection form at finitely many levels. We show that in each residue class these levels form an interval, and that once $hd\ge c_1+c_2$ the total of $I_h$ is exactly the mass, so that every sum missing from $hA_0$ raises the total above it.

For a set $E$ of integers containing $0$ and an integer $n\ge0$ put $D_E(n)=n-\max(E\cap[0,n])$, so that $D_h=D_{hA_0}$. For every $A_0$ and $h\ge0$ we have $hA_0\subseteq I_h$ and $0,hd\in hA_0$. Put $F_h=I_h\setminus hA_0$, and call $h$ a \emph{failing level} if $F_h\neq\varnothing$. By Theorem~\ref{thm:cm}, $A_0$ is CM if and only if it has no failing level.

\begin{proposition}[Failing levels]\label{prop:failing-interval}
Fix $1\le s\le d-1$. The levels $h$ at which some integer $kd-s$ lies in $F_h$ form an interval of integers, possibly empty. Consequently $A_0$ has finitely many failing levels.
\end{proposition}

\begin{proof}
\emph{The corners.} Write a target of class $-s$ as $kd-s$ and use coordinates $(k,\ell)$ with $\ell=h-k$. Let $\mathcal C_s$ be the set of $(k,\ell)\in\N^2$ with $kd-s\in(k+\ell)A_0$, and let $\mathcal C^I_s$ be the set with $kd-s\in I_{k+\ell}$. By Lemma~\ref{lem:residue} with $m=d$, $\mathcal C^I_s=[u_s,\infty)\times[v_s,\infty)$. The set $\mathcal C_s\subseteq\mathcal C^I_s$ is closed under adding $(1,0)$ (add a part $d$) and $(0,1)$ (add a part $0$), so it is the union of the quadrants above its finitely many minimal points (Dickson's lemma) $(\kappa_1,\mu_1),\dots,(\kappa_t,\mu_t)$, ordered so that $\kappa_1<\dots<\kappa_t$ and $\mu_1>\dots>\mu_t$. We claim $\kappa_1=u_s$ and $\mu_t=v_s$. The Ap\'ery element $u_sd-s$ of $S_1$ is a sum of some number $j$ of positive parts, so $u_sd-s\in jA_0$; as all parts are at most $d$ and $u_sd-s>(u_s-1)d$, we have $j\ge u_s$, so $(u_s,j-u_s)\in\mathcal C_s$. Likewise the Ap\'ery element $v_sd+s$ of $S_2$ is a sum of $j'$ parts $d-a$, $a\in A_0$, so $j'd-(v_sd+s)=(j'-v_s)d-s$ lies in $j'A_0$, and $(j'-v_s,v_s)\in\mathcal C_s$. Together with $\mathcal C_s\subseteq\mathcal C^I_s$ this proves the claim.

\emph{The failing levels of one class.} If $t=1$, then $\mathcal C_s=\mathcal C^I_s$ and no level fails in this class. If $t\ge2$, then $\mathcal C^I_s\setminus\mathcal C_s$ is the disjoint union of the rectangles $[\kappa_j,\kappa_{j+1}-1]\times[v_s,\mu_j-1]$, $1\le j<t$, because for $\kappa_j\le k<\kappa_{j+1}$ the least $\ell$ with $(k,\ell)\in\mathcal C_s$ is $\mu_j$. The levels $k+\ell$ met by rectangle $j$ form the interval $J_j=[\kappa_j+v_s,\ \kappa_{j+1}+\mu_j-2]$, which is nonempty. For $j\le t-2$, strict descent gives $\mu_j\ge\mu_{j+1}+1\ge v_s+2$, so the upper end of $J_j$ is at least $\kappa_{j+1}+v_s$, the lower end of $J_{j+1}$. Hence $J_1\cup\dots\cup J_{t-1}$ is the interval $[u_s+v_s,\ \max_{j<t}(\kappa_{j+1}+\mu_j-2)]$. Multiples of $d$ never fail, so the failing levels are the finite union of these intervals over $s$.
\end{proof}

In other coordinates, Proposition~\ref{prop:failing-interval} and the formula for $\tau$ below are \cite[Lemma~2.5]{CompanionRegularity}; we keep the short proof.

Let $\tau$ be the least $h_0\ge0$ such that $hA_0=I_h$ for every $h\ge h_0$. Thus $\tau=0$ for CM sets, and otherwise $\tau-1$ is the last failing level; by the proof above, $\tau-1=\max_s\max_{j<t_s}(\kappa_{j+1}+\mu_j-2)$ over the classes $s$ with $t_s\ge2$ corners, an explicit finite computation.

\begin{proposition}[Overshoot]\label{prop:overshoot}
Let $h\ge0$ satisfy $hd\ge c_1+c_2$. Then $\sum_{n=0}^{hd}D_{I_h}(n)=\M(S_1)+\M(S_2)$. Consequently $\Tot(h)\ge\M(S_1)+\M(S_2)$, with equality if and only if $h$ is not a failing level. In particular, if $A_0$ is not CM and $(\tau-1)d\ge c_1+c_2$, then $\Tot(\tau-1)>\M(S_1)+\M(S_2)$.
\end{proposition}

\begin{proof}
Since $hd-c_1\ge c_2$, the integer $c_1$ lies in $I_h$, and so does every integer of $[c_1,hd-c_2]$. An integer $n<c_1$ has $hd-n>c_2$, so $n\in I_h$ if and only if $n\in S_1$; an integer $n>hd-c_2$ has $n>c_1$, so $n\in I_h$ if and only if $hd-n\in S_2$. Thus the integers of $[0,hd]$ outside $I_h$ are the gaps of $S_1$, all below $c_1$, and the reflections $hd-x$ of the gaps $x$ of $S_2$, all above $hd-c_2$, and the element $c_1$ of $I_h$ separates the two kinds. A maximal run of $L$ gaps of $S_1$ contributes $1+\dots+L$ to the total of $I_h$ and to $\M(S_1)$. For a reflected run, the undershoot at $hd-x$ is the distance from $x$ up to the next element of $S_2$, so by Lemma~\ref{lem:elementary}(iv) the run contributes $L(L+1)/2$, as in $\M(S_2)$. This proves the first identity (if a conductor is $0$, the corresponding set of gaps is empty). Removing from a set containing $0$ and $hd$ an element $x$ with $0<x<hd$ makes its undershoot at $x$ positive and decreases no other undershoot (Lemma~\ref{lem:elementary}(iii)); since $hA_0=I_h\setminus F_h$ and $0,hd\notin F_h$, the total of $hA_0$ exceeds that of $I_h$ exactly when $F_h\neq\varnothing$.
\end{proof}

So once the range is long enough, the intersection form carries exactly the mass, and each sum that $hA_0$ lacks adds to the total. For a CM set Proposition~\ref{prop:overshoot} recovers $\Tot(h)=\M(S_1)+\M(S_2)$ whenever $hd\ge c_1+c_2$. For four-element sets the hypothesis on $\tau$ always holds: by \cite[Theorem~A]{CompanionRegularity}, every non-CM four-element set satisfies
\[
(\tau-1)d\ge c_1+c_2+d ,
\]
with equality exactly for $\{0,b-1,b+1,2b\}$ with $b\ge2$ even.

\begin{corollary}\label{cor:four-overshoot}
Every four-element set $A_0$ that is not CM satisfies $\Tot(\tau-1)>\M(S_1)+\M(S_2)$. In particular, for four-element sets, $\Tot(h)\le\M(S_1)+\M(S_2)$ for all $h$ implies that $A_0$ is CM.
\end{corollary}

\begin{proof}
Combine \cite[Theorem~A]{CompanionRegularity} with Proposition~\ref{prop:overshoot}.
\end{proof}

With five elements the conclusion can fail.

\begin{example}\label{ex:five-element}
Let $A_0=\{0,5,9,11,20\}$, the non-Cohen--Macaulay set of \cite[Example~4.3]{BermejoGimenez2006} (see also \cite[Example~2.13(2)]{GimenezGonzalezSanchez2023}). Then $c_1=18$, $c_2=44$ and $\M(S_1)+\M(S_2)=74$. The only failing level is $h=3$, so $A_0$ is not CM and $\tau=4$, and $\Tot(h)=53,69,74$ for $h=1,2,3$. For $h\ge4$ we have $hA_0=I_h$ and $20h\ge80>c_1+c_2=62$, so $\Tot(h)=74$ by Proposition~\ref{prop:overshoot}. Hence $\Tot(h)\le\M(S_1)+\M(S_2)$ for every $h$: the set never overshoots its mass. At the failing level $h=3$ the hypothesis $hd\ge c_1+c_2$ fails ($60<62$), and the total equals the mass although $F_3=\{24\}$: the total of $I_3$ is $73$, and the correction of Proposition~\ref{prop:correction} at $z=24$, with $z^-=23$ and $z^+=25$, adds $1$.
\end{example}

\section{Corrections and windows}\label{sec:corrections}

This section treats arbitrary sets and finite windows of targets. The undershoot of an arbitrary set at any level is that of $I_h$ plus an explicit correction, and the undershoot over a window has a closed form in the corners of Section~\ref{sec:noncm}. With $D_E$ and $F_h$ as in Section~\ref{sec:noncm} put
\[
o_2(y)=\min\{j\ge0:y+j\in S_2\}.
\]
An \emph{equality level} is an $h$ with $hA_0=I_h$: every $h$ for Cohen--Macaulay sets, every $h\ge\tau$ otherwise.

\subsection{Undershoot as overshoot}\label{subsec:overshoot-form}

\begin{proposition}[Head and reflected tail]\label{prop:overshoot-form}
Let $A_0$ be arbitrary and $h\ge0$.
\begin{enumerate}
\item[(a)] $D_{I_h}(n)=D_{S_1}(n)$ for $0\le n\le hd-c_2$.
\item[(b)] If $x\in I_h$ and $x\ge c_1$, then $D_{I_h}(n)=o_2(hd-n)$ for $x\le n\le hd$. One may take $x=c_1$ when $hd-c_1\in S_2$, in particular when $c_1=0$ or $hd\ge c_1+c_2$.
\item[(c)] If $hd\ge c_1>0$ and $N_h=\max(I_h\cap[0,c_1-1])$, then for $c_1\le n\le hd$, $D_{I_h}(n)=o_2(hd-n)$ if $o_2(hd-n)\le n-c_1$ and $D_{I_h}(n)=n-N_h$ otherwise.
\item[(d)] If $hd\ge c_1+c_2-1$, (a) together with (b) or (c) determines $D_{I_h}$ on $[0,hd]$.
\end{enumerate}
At an equality level $D_h=D_{I_h}$.
\end{proposition}

\begin{proof}
(a) For $n'\le n\le hd-c_2$, $hd-n'\ge c_2$, so $n'\in I_h$ if and only if $n'\in S_1$. (b), (c) For $c_1\le n'\le n$, $n'\in S_1$, so $n'\in I_h$ if and only if $hd-n'\in S_2$. With $j=o_2(hd-n)$: if $j\le n-c_1$, then $n-j$ is the largest element of $I_h$ in $[c_1,n]$; otherwise there is none and the predecessor is $N_h$ (and $0\in I_h$). If $x\in I_h$, $x\ge c_1$ and $n\ge x$, then $j\le n-x$, since $(hd-n)+(n-x)=hd-x\in S_2$. If $c_1=0$, $hd-c_1=hd\in S_2$. (d) When $hd\ge c_1+c_2-1$ we have $hd-c_2\ge c_1-1$, so the range $[0,hd-c_2]$ of (a) and the range $[c_1,hd]$ of (b) and (c) together cover $[0,hd]$.
\end{proof}

So above the anchor the undershoot is the overshoot of the reflected target into $S_2$; summed over $[0,hd]$ this is Proposition~\ref{prop:overshoot}. For an arbitrary set the actual undershoot differs from that of $I_h$ by an explicit correction.

\begin{proposition}[Correction]\label{prop:correction}
For $z\in F_h$ let $z^+=\min(I_h\cap(z,hd])$ and $z^-=\max(hA_0\cap[0,z))$. For every $W\subseteq[0,hd]$,
\[
\sum_{n\in W}D_h(n)=\sum_{n\in W}D_{I_h}(n)+\Psi_h(W),\qquad \Psi_h(W)=\sum_{z\in F_h}\#\bigl(W\cap[z,z^+)\bigr)(z-z^-).
\]
So $\Psi_h\ge0$, and $\Psi_h([0,hd])>0$ if and only if $h$ is a failing level.
\end{proposition}

\begin{proof}
Both $z^\pm$ exist as $0,hd\in hA_0$. For a target $n$ let $z=\max(I_h\cap[0,n])$. If $z\notin F_h$ it is also the predecessor in $hA_0\subseteq I_h$. If $z\in F_h$, then $n\in[z,z^+)$, $hA_0$ has no element in $[z,n]$, and the predecessor is $z^-$, adding $z-z^-$. The cells $[z,z^+)$ are disjoint.
\end{proof}

For weighted targets replace the cardinality by the weight of $W\cap[z,z^+)$.

\subsection{Window formulas}\label{subsec:overlap}

For arbitrary sets the correction on a window has a finite closed form in the corners of Proposition~\ref{prop:failing-interval}. Let $C_s$ be the set of minimal points of $\mathcal C_s$, so that $kd-s\in hA_0$ if and only if $\kappa\le k\le h-\mu$ for some $(\kappa,\mu)\in C_s$. These corners are, in other coordinates, the strict records of the least number of parts in the class of $-s$ that Chen and Yang use, and for a single position Lemma~\ref{lem:hinge} is the positive-part count behind their exact formula for $|hA_0|$ \cite[Lemma~2, Thm.~1]{ChenYang2026}; the window clipping, the union over consecutive positions and the passage to the undershoot and its correction are what Theorem~\ref{thm:hinge} adds. For $1\le s\le t<d$ let $C^B_{s,t}$ be the minimal elements of $\bigcup_{r=s}^tC_r$ and $C^I_{s,t}$ those of $\{(u_r,v_r):s\le r\le t\}$, listed with $\kappa_1<\dots<\kappa_J$ and $\mu_1>\dots>\mu_J$. Each corner $(\kappa_j,\mu_j)$ makes a position present in the blocks $k\in[\kappa_j,h-\mu_j]$, and the window keeps the blocks in $[\ell_0,q_0]$. The function $H_C$ below counts the blocks of $[\ell_0,q_0]$ covered by at least one of these intervals: it adds their clipped lengths and subtracts the overlap of each interval with the one before it, which is all of its overlap with the earlier ones, since both endpoint sequences increase. With $[x]_+=\max(x,0)$ put
\begin{multline*}
H_C(h;\ell_0,q_0)=\sum_{j=1}^{J}\bigl[\min(q_0,h-\mu_j)-\max(\ell_0,\kappa_j)+1\bigr]_+\\
-\sum_{j=2}^{J}\bigl[\min(q_0,h-\mu_{j-1})-\max(\ell_0,\kappa_j)+1\bigr]_+ .
\end{multline*}

\begin{lemma}\label{lem:hinge}
$H_C(h;\ell_0,q_0)$ is the number of integers in $[\ell_0,q_0]\cap\bigcup_j[\kappa_j,h-\mu_j]$.
\end{lemma}

\begin{proof}
The intervals $I_j=[\max(\ell_0,\kappa_j),\min(q_0,h-\mu_j)]$ have nondecreasing endpoints, so a point of $I_j\cap I_i$, $i<j-1$, lies in $I_{j-1}$, and $I_j\cap(I_1\cup\dots\cup I_{j-1})=I_j\cap I_{j-1}$, of size given by the $j$-th subtracted term (also for empty intervals).
\end{proof}

\begin{theorem}[Hinge formulas]\label{thm:hinge}
Let $T_-\le T_+$ be integers in $[0,hd]$, and put $\ell_s=\max(1,\lceil(T_-+s)/d\rceil)$ and $q_s=\min(h,\lfloor(T_++s)/d\rfloor)$. Then
\begin{gather*}
\sum_{T=T_-}^{T_+}D_h(T)=\sum_{1\le s\le t<d}\Bigl([q_s-\ell_s+1]_+-H_{C^B_{s,t}}(h;\ell_s,q_s)\Bigr),\\
\Psi_h([T_-,T_+])=\sum_{1\le s\le t<d}\bigl(H_{C^I_{s,t}}(h;\ell_s,q_s)-H_{C^B_{s,t}}(h;\ell_s,q_s)\bigr).
\end{gather*}
For $0<c<d$, $c\notin A_0$, with corners $C'$ for $A_0\cup\{c\}$, the total undershoot of $A_0\cup\{c\}$ at level $h$ on $[T_-,T_+]$ is smaller than that of $A_0$ at level $h+1$ if and only if $\sum_{s\le t}H_{C'^B_{s,t}}(h;\ell_s,q_s)>\sum_{s\le t}H_{C^B_{s,t}}(h+1;\ell_s,q_s)$, and they tie exactly when the sums are equal.
\end{theorem}

\begin{proof}
Multiples $kd$, $k\le h$, lie in $hA_0$, so as in Proposition~\ref{prop:counting}(a), $D_h(kd-s)$ counts the $t\ge s$ such that no position of $[s,t]$ is present in block $k$. The target $kd-s$ lies in $[T_-,T_+]$ exactly for $\ell_s\le k\le q_s$. Position $r$ is present in block $k$ exactly for $k\in[\kappa,h-\mu]$, $(\kappa,\mu)\in C_r$, and dominated corners give smaller intervals; Lemma~\ref{lem:hinge} counts the blocks in which some position of $[s,t]$ is present. For $I_h$ the corners are the $(u_r,v_r)$ (Lemma~\ref{lem:residue}). For the comparison, $T_+\le hd$ and $s<d$ give $\lfloor(T_++s)/d\rfloor\le h$, so $\ell_s,q_s$ are the same at levels $h$ and $h+1$ and the first terms cancel.
\end{proof}

The corners depend neither on $h$ nor on the window, so finitely many integers per residue class give the undershoot over every window at every level.
For $[T_-,T_+]=[0,hd]$ every term is $[h+1-\kappa-\mu]_+$ with fixed corners, so $\Psi_h([0,hd])$ is piecewise affine in $h$ and vanishes for $h\ge\tau$.

\section{Open problems}\label{sec:open}

\begin{enumerate}
\item \emph{Similar ordering through continued fractions.} Theorem~\ref{thm:similar} characterizes the similarly ordered Cohen--Macaulay four-element sets by the two Ap\'ery sets. Give a characterization by a short list of congruence or continued-fraction conditions on $(a_1,a_2,d)$, for instance through the Euclidean-type algorithm that computes the Ap\'ery sets of the three-generated semigroups $S_1$ and $S_2$ \cite{Rodseth1978}, and extend the diagram description and the test to five or more elements.
\item \emph{Other families.} The spectrum of $\{0,a,b,2b\}$ is governed by the three-distance theorem (Section~\ref{sec:family}). Compute the spectrum at the modulus $d$ of other complete-intersection families, such as $\{0,a,b,kb\}$ with $0<a<b$ coprime and $k\ge3$; for $k\ge3$, $b$ is not a modulus, since $b\notin S_2$.
\item \emph{Overshoot for more elements.} Every non-Cohen--Macaulay four-element set overshoots its mass at some level, but $\{0,5,9,11,20\}$ does not (Section~\ref{sec:noncm}). Which non-Cohen--Macaulay sets with five or more elements overshoot their mass? Proposition~\ref{prop:overshoot} and Theorem~\ref{thm:hinge} decide it for any given set by a finite computation; the question asks for a structural criterion. Hoa and Tien \cite[Thm.~4.4]{HoaTien2026} describe the points of $I_h$ outside $hA_0$, as the first local cohomology of $\Bbbk[\mathcal S]$, through Ap\'ery sets, which may be a starting point.
\end{enumerate}

\section{Computation and reproducibility}\label{sec:computation}

Every computational statement in this paper is reproduced by the accompanying code \cite{code2026}, which uses only the Python standard library and exact integer and rational arithmetic. Three independent engines serve as oracles: direct iterated sumsets, the homogenized Ap\'ery set (from which $hA_0$ and the undershoot are rebuilt and compared with the first engine), and the intersection form of Theorem~\ref{thm:cm} computed from the two numerical semigroups. The regularity is computed from \eqref{eq:reg-cm}, not by a computer-algebra system; the generators of the ideal in Lemma~\ref{lem:family-ci} are checked separately by Gr\"obner bases in SymPy \cite{SymPy2017} for $b\le12$, over $\mathbb Q$, $\mathbb F_2$ and $\mathbb F_3$; and the inequality of \cite[Theorem~A]{CompanionRegularity} is checked only as finite evidence. Every formula of the paper was checked on the ranges in Table~\ref{tab:checks}: against direct sumsets and the definitions where the table gives a direct range (for the family $b\le25$), and beyond that against the other exact formulas of the paper; the outputs are deterministic, and a manifest of their SHA-256 digests lets a reader confirm a byte-identical rerun. The full suite runs in under a minute on one core.

\begin{table}[ht]
\caption{Exact checks. ``Sets with $d\le D$'' means all $A_0=\{0<a_1<\dots<d\}$ with $\gcd(A_0)=1$ and the stated number of elements.}\label{tab:checks}
\begin{tabular}{@{}p{0.44\textwidth}p{0.52\textwidth}@{}}
\toprule
Results & Range checked \\
\midrule
Sections~\ref{sec:intro}, \ref{sec:sumsets} and~\ref{sec:similar} & sets with three and four elements, $d\le20$; five elements, $d\le12$; census $d\le30$ \\
Section~\ref{sec:family} & $b\le25$ by direct sumsets, $b\le60$ by formula, three distances $b\le80$; ideal of Lemma~\ref{lem:family-ci} by Gr\"obner bases $b\le12$ \\
Section~\ref{sec:noncm} & four elements $d\le23$, five elements $d\le14$ \\
Section~\ref{sec:corrections} & sets with $d\le11$ (five elements $d\le8$), levels $h\le5$; window formulas $d\le10$, $h\le4$ \\
\botrule
\end{tabular}
\end{table}

No result of the paper depends on computation; the checks confirm the formulas on the stated ranges.

\subsection*{Use of generative AI}\label{sec:ai}
The authors used Claude Opus (Anthropic) to search the literature, to draft proofs and text under the authors' direction, to write the verification code in the accompanying archive, and to write the Ti\emph{k}Z code of Figure~\ref{fig:diagrams} from the computed Ap\'ery sets. The authors checked all output and take full responsibility for the content.

\backmatter

\section*{Statements and Declarations}

\begin{itemize}
\item \textbf{Funding.} No funding was received for conducting this study.
\item \textbf{Competing interests.} The authors have no competing interests to declare that are relevant to the content of this article.
\item \textbf{Ethics approval and consent to participate.} Not applicable.
\item \textbf{Consent for publication.} Not applicable.
\item \textbf{Data availability.} All data generated by the computations are archived on Zenodo at \url{https://doi.org/10.5281/zenodo.22952563} \cite{code2026}.
\item \textbf{Materials availability.} Not applicable.
\item \textbf{Code availability.} The code is archived with the data at the same DOI, from the tagged release of the source repository.
\end{itemize}

\bibliography{refs}

\end{document}